\documentclass[12pt]{amsart}
\usepackage[english]{certus}
\usepackage{amsrefs, setspace, a4wide}

\title{Quadratic Modules, $C^*$-Algebras, and Free Convexity}
\author{Vadim Alekseev}
\address{V.A., Institut f\"ur Geometrie, TU Dresden, 01062 Dresden, Germany}
\email{vadim.alekseev@tu-dresden.de}

\author{Tim Netzer}
\address{T.N., Institut f\"ur Geometrie, TU Dresden, 01062 Dresden, Germany}
\email{tim.netzer@tu-dresden.de}

\author{Andreas Thom}
\address{A.T., Institut f\"ur Geometrie, TU Dresden, 01062 Dresden, Germany}
\email{andreas.thom@tu-dresden.de}

\begin{document}

\onehalfspace
 \maketitle
 
 \begin{abstract}
Given a quadratic module, we construct its universal $C^*$-algebra, and then use methods and notions from the theory of $C^*$-algebras to study the quadratic module.  We define residually finite-dimensional quadratic modules, and characterize them in various ways, in particular via a Positivstellensatz. We give unified proofs for several existing strong Positivstellens\"atze, and prove some new ones. Our approach also leads naturally to interesting new examples in free convexity.  We show that the usual notion of a free convex hull is not able to detect residual finite-dimensionality. We thus propose a new notion of free convexity, which is coordinate-free. We characterize semialgebraicity of free convex hulls of semialgebraic sets, and show that they are not always semialgebraic, even at scalar level. This also shows that the membership problem for quadratic modules has a negative answer in the non-commutative setup.
\end{abstract}

\section{Introduction}

Quadratic modules are well-studied objects in real algebra. They are generalizations of the cone of sums of squares, and play the role in Positivstellens\"atzen that ideals play in Nullstellens\"atzen. The commutative theory of quadratic modules is quite well-understood (see \cite{bcr,ma,pd} or \cite{sch} for a survey). Interest in the non-commutative theory is much more recent (see \cite{schm} for a survey; more references to non-commutative Positivstellens\"atze can be found throughout this article).  
Also a quite new development with many recent results is free convexity (see for example \cite{eff,he1,he3,he4}). Instead of looking at convex  sets in $\mb R^d$ only, one considers sets of matrix-tuples of all sizes simultaneously. A suitable notion of convexity then relates the different matrix levels. All these notions are well-motivated by applications in such diverse areas as quantum physics, linear systems engineering, free probability and semidefinite optimization (see \cite{he1} for more information). 

In this paper, our contribution is the following.  In Section \ref{cstar} we consider the $C^*$-algebra that one can canonically assign to an (archimedean) quadratic module (this was also done in \cite{ci}). With this construction, some of the most important methods from operator algebra pass to quadratic modules, as we explain. After assembling the necessary techniques, we define residually finite-dimensional (r.f.d.) quadratic modules in Section \ref{rfdsec}. This notion exists for $C^*$-algebras, and has interesting characterizations when formulated for quadratic modules. It corresponds to a Positivstellensatz with positivity at finite-dimensional representations. In this context we give alternative and uniform proofs for the strong Positivstellens\"atze from \cite{he2, he5}, and prove the same results for more classes of examples. In Section \ref{convex} we investigate free convexity. The results from Section \ref{rfdsec} indicate that the present notions in free convexity are not always optimal. The coordinate-based approach is for example not always able to detect the property r.f.d.\ of a quadratic module.  We thus suggest a new and coordinate-free approach towards free convexity. We characterize free convex hulls of semialgebraic sets to be semialgebraic at each matrix level. We produce examples showing that both matrix- and operator-convex hulls of free semialgebraic sets can fail to be semialgebraic, even at scalar level. This also shows that for a finitely generated quadratic module in a free algebra, the intersection with a finite-dimensional subspace can fail to be semialgebraic. So the membership problem has a negative answer in the non-commutative setup (see \cite{aug} for partial positive answers in the commutative case).

\section{Universal $C^*$-algebras of quadratic modules}\label{cstar}

In this section we define the most important notions of the paper, and assemble important techniques for later use. 
\begin{Def}
A quadratic module $(\mc A, \mc Q)$ is a pair consisting of a unital complex $*$-algebra and a subset $\mathcal{Q}\subseteq \mathcal A^h$, where $\mathcal A^h$ denotes the $\mathbb R$-subspace of hermitian elements of $\mathcal A$, such that $1\in \mathcal{Q}$ and $$a,b\in \mathcal{Q}, c\in \mathcal A\quad \Rightarrow\quad  c^*(a+b)c \in \mathcal{Q}.$$
\end{Def}

\begin{Def}
The quadratic module $(\mc A,\mc Q)$ is archimedean if $\ell-a^*a\in \mathcal{Q}$ for any $a\in\mathcal A$ and large enough $\ell$, or equivalently $\ell-a\in \mathcal{Q}$ for any $a\in\mathcal A^h$ and large enough $\ell$. It is also enough to require this for generators of $\mathcal A$ only (see \cite{ci} for technical details). 
\end{Def}

It's clear from the definitions that any quadratic module contains sums of squares 
\[
 \Sigma^2\mc A = \left\{\sum_{i=1}^n a_i^*a_i\,|n\in\mb N,\, a_i\in\mc A\right\}.
\]
In the sequel, if we don't specify a quadratic module in a $\ast$-algebra $\mc A$, we always assume that it comes with the smallest quadratic module $\Sigma^2\mc A$. Notice that for a \Cs-algebra $\mc A$, $\Sigma^2 \mc A$ is just the set of positive elements of $\mc A$ (for an introduction to $C^*$-algebras consult for example \cite{arv1}).

\begin{Def}
 Let $(\mc A,\mc Q)$ and $(\mc B, \mc R)$ be two quadratic modules. Their tensor product
 \[
  (\mc A,\mc Q)\otimes (\mc B,\mc R)
 \]
is defined to be the smallest quadratic module in $\mc A\otimes\mc B$ containing the set $\{q\otimes r\,|\, q\in\mc Q,\,r\in \mc R\}$.
\end{Def}

\begin{Ex}(1)
 If $(\mc B,\mc R) = \mb M_n(\mb C)$, then
 \[
  (\mc A,\mc Q)\otimes \mb M_n(\mb C) = \left(\mb M_n(\mc A), \left\{ \sum_{\mbox{\tiny finite}} \left( a_i^*qa_j\right)_{ij}\mid a_1,\ldots, a_n\in \mathcal A, q\in \mathcal{Q}\right\}\right).
 \]
For sake of brevity we denote this quadratic module also by $\mc Q\otimes \mb M_n$. It is not hard to check that $\mc Q\otimes \mb M_n$ is archimedean, if $\mc Q$ was.

(2) Let $\mb R\subseteq {\bf R}$ be an extension of real closed fields. Then $\mb C\subseteq {\bf C}={\bf R} [i]$ is an extension of algebraically closed fields, and $\bf C$ is a unital $*$-algebra over $\mb C$ with ${\bf C}^h=\bf R$. For $\mc Q=\Sigma^2{\bf C}={\bf R}_+$ we find $\mc Q\otimes\mb M_n=\mb M_n({\bf C})_+,$ the set of positive semidefinite hermitian matrices over $\bf C$. Note that semidefiniteness has the same characterizations over ${\bf C}$ as over $\mb C$, for example by Tarski's Transfer Principle (see for example \cite{pd} for details on real closed fields and their model-theoretic properties).
\end{Ex}

\begin{Def}
Let $(\mathcal A, \mathcal Q)$ and $(\mathcal B, \mathcal R)$ be two quadratic modules. 

\noindent
(i) A unital completely positive morphism (u.c.p. morphism)
\[
 \varrho\colon (\mathcal A, \mathcal Q) \to (\mathcal B, \mathcal R)
\]
is a unital $*$-linear map $\varrho\colon \mathcal A\to\mathcal B$ such that for every $n\in \mb N$,
\[
(\varrho\otimes\id)(\mathcal Q\otimes \mb M_n)\subseteq \mathcal R \otimes \mb M_n.
\]

\noindent
(ii) A homomorphism
\[
 \pi\colon (\mathcal A, \mathcal Q) \to (\mathcal B, \mathcal R)
\]
is a unital $\ast$-homomorphism $\pi\colon \mathcal A\to\mathcal B$ such that $\pi(\mathcal Q)\subseteq \mathcal R$. 

\noindent
(iii) A representation of a quadratic module $(\mc A,\mc Q)$ is a homomorphism
\[
 \pi\colon (\mc A,\mc Q) \to \mb B(\mc H),
\]
where $\mc H$ is a Hilbert space.

\end{Def}

Notice that any homomorphism between quadratic modules is obviously a u.c.p. morphism.
Of course, a representation of $(\mc A,\mc Q)$ is just a $*$-representation $\pi\colon \mc A \to \mb B(\mc H)$ fulfilling $\pi(q)\geqslant 0$ for all $q\in \mathcal{Q}$. We will see that representations of a quadratic module are in one-to-one correspondence with the representations of its universal \Cs-algebra. 

Although the collection of all representations of a quadratic module is in general not a set, this problem can be avoided by appropriately resticting the cardinality of the target Hilbert space, for instance, bounding it to the cardinality of the universal representation of the corresponding \Cs-algebra. We will tacitly do this and denote by $\Rep(\mc A,\mc Q)$ the set of all representations of a quadratic module $(\mc A,\mc Q)$ and by $\Rep_{\mathrm{fd}}(\mc A,\mc Q)$ the subset of finite-dimensional representations.

\begin{Def}
Let $(\mc A,\mathcal Q)$ be an archimedean quadratic module. We equip the algebra $\mathcal A$ with the seminorm
\[
 \norm{a}_{\mathcal Q} = \sup_{\pi \in \Rep(\mc A,\mathcal Q)} \norm{\pi(a)}.
\]
The supremum is finite, because
\[
 \norm{\pi(a)}^2 = \norm{\pi(a^*a)} \leqslant \ell,\quad \pi\in\Rep(\mc A,\mc Q),
\]
if $\ell - a^*a \in \mathcal Q$. In fact we have $$\Vert a\Vert_{\mathcal Q}=\inf\left\{\ell \mid \ell^2-a^*a\in \mathcal Q\right\},$$ which follows by separation from the archimedean cone $\mathcal Q$, and the GNS construction.  

As in \cite{ci}, the (separated) completion of $\mathcal A$ with respect to $\norm{\cdot}_{\mc Q}$ is denoted by $C^*(\mathcal A,\mathcal Q)$ and called the universal $C^*$-algebra of $(\mc A,\mathcal Q)$. We denote by $\iota\colon \mathcal A \to C^*(\mathcal A,\mathcal Q)$ the canonical map with dense image.
\end{Def}

The name ``universal'' comes from the following fact:
\begin{Prop}\label{universal}
Let $(\mathcal A,\mathcal Q)$  be an archimedean quadratic module, $\iota\colon \mathcal A \to C^*(\mathcal A,\mathcal Q)$ the canonical map, and $\mc B$ a $C^*$-algebra.

(i) $\iota$ is a homomorphism of quadratic modules and respects the (semi)-norm. 


(ii) There is a one-to-one correspondence between u.c.p morphisms $\varrho\colon(\mc A,\mc Q)\rightarrow \mc B$ and u.c.p. morphisms $\overline\varrho\colon C^*(\mc A,\mc Q)\rightarrow \mc B$. The correspondence is given by the formula
 $\varrho = \overline\varrho \circ \iota.$ This correspondence maps homomorphisms to homomorphisms.

\end{Prop}
\begin{proof}($i$) For any $\ell>\Vert q\Vert_{\mathcal Q}$ we have $\ell^2-q^*q\in\mathcal Q$. Then $$\ell-q=\frac{1}{2\ell}\left((\ell^2-q^*q)+(\ell-q)^*(\ell-q)\right) \in\mathcal Q.$$ In particular $2\ell-q\in\mathcal Q$ and thus $$\ell^2-(\ell-q)^*(\ell-q)=\frac{1}{2\ell}\left((2\ell-q)^*q(2\ell-q)+ q^*(2\ell-q)q\right)\in\mathcal Q.$$ This proves $\Vert \ell-q\Vert_{\mathcal Q}\leq\ell$, whenever $\ell>\Vert q\Vert_{\mathcal Q}$. The same then holds in $C^*(\mathcal A,\mathcal Q)$, and this is well known to imply positivity in a $C^*$-algebra.


For ($ii$) first note that $\iota$ is u.c.p., and any u.c.p. map on $C^*(\mc A,\mc Q)$ is continuous. Further, any u.c.p. map $\varrho$ on $(\mc A,\mc Q)$ factors through $C^*(\mc A,\mc Q)$. Indeed if $\Vert a\Vert_{\mc Q}=0$, then $\epsilon-a^*a\in \mc Q$ for all $\epsilon>0$. Then $$\left(\begin{array}{cc}1 & a \\a^* & \epsilon\end{array}\right) = \left(\begin{array}{c}1 \\a^*\end{array}\right)  \left(\begin{array}{cc}1 & a\end{array}\right) + \left(\begin{array}{c}0 \\1\end{array}\right)(\epsilon -a^*a)\left(\begin{array}{cc}0 & 1\end{array}\right)\in \mc Q\otimes \mb M_2$$  and thus $$\left(\begin{array}{cc}1 & \varrho(a) \\\varrho(a)^* & 0\end{array}\right)\geqslant 0,$$ which implies $\varrho(a)=0$. Finally,  $(\iota\otimes {\rm id})(\mc Q\otimes \mb M_n)$ is dense in the positive elements of $C^*(\mc A,\mc Q)\otimes \mb M_n(\mb C),$ so $\overline\varrho$ is u.c.p.
\end{proof}

We now formulate some important techniques from operator algebra in the context of quadratic modules (see \cite{pau} for the corresponding results for $C^*$-algebras).  

\begin{Prop}[Stinespring's Dilation Theorem]\label{stinespring}
Let $(\mc A,\mathcal{Q})$ be an archimedean quadratic module and let
\[
 \varrho\colon (\mc A,\mc Q) \to \mb B(\mc H_\varrho)
\]
be a u.c.p. morphism. Then there is a representation $\pi \colon (\mc A,\mc Q) \to \mb B(\mc H_\pi)$ and an isometry $\gamma\colon \mathcal H_\varrho\hookrightarrow \mathcal H_{\pi}$ such that $\varrho= {}^{\gamma} \pi$, where 
\[
{}^\gamma \pi(a)= (\gamma^*\circ\pi(a)\circ\gamma).
\]
 \end{Prop}
 \begin{proof}
For $C^*$-algebras, this is precisely the statement of Stinespring's Dilation Theorem (\cite{pau}, Theorem 4.1). The version for quadratic modules is immediate from Proposition \ref{universal}. \end{proof}
 
 Now let $\mathcal V\subseteq \mathcal A$ be a unital $*$-subspace. For $n\in\mathbb N$ we equip the $*$-space $\mc V\otimes {\rm Mat}_n(\mb C)$ with the convex cone $$(\mc Q\otimes \mb M_n)_{\vert \mc V}:=\left(\mathcal Q\otimes\mb M_n\right)\cap\left(\mc V\otimes  {\mb M}_n(\mb C)\right).$$ A unital $*$-linear mapping $\varrho\colon \mathcal V\rightarrow \mathcal B(\mathcal H)$ is again called  u.c.p. if all mappings $\varrho\otimes{\rm id}$ map these cones to positive elements.

\begin{Prop}[Arveson's Extension Theorem]\label{extend} Let $(\mc A,\mathcal Q)$ be an archimedean quadratic module and  $\mathcal V\subseteq \mathcal A$  a unital $*$-subspace. Then any u.c.p.\ map $\varrho\colon \mathcal V\rightarrow \mathcal B(\mathcal H)$ extends to a u.c.p.\ map $\varrho\colon \mathcal A\rightarrow \mathcal B(\mathcal H).$
\end{Prop}
\begin{proof} Any u.c.p.\ map $\varrho\colon \mc V\rightarrow\mb B(\mc H)$ factors through $\iota(\mc V)$, by the same argument as in Proposition \ref{universal}. We show that the resulting map $\overline\varrho\colon \iota(\mc V)\rightarrow \mb B(\mc H)$ is u.c.p., and the result then  clearly follows from the standard version of Arveson's Extension Theorem \cite[Theorem 7.5]{pau}.

It is not hard to check that $C^*(\mc A\otimes \mb M_n(\mb C),\mc Q\otimes \mb M_n)=C^*(\mc A,\mc Q)\otimes\mb M_n(\mb C)$ holds. Every $*$-linear functional $\varphi\colon \mc A\otimes\mb M_n(\mb C)\rightarrow\mb C$ which is nonnegative on $\mc Q\otimes\mb M_n$ is automatically u.c.p.,  and thus extends to a u.c.p. functional on $C^*(\mc A,\mc Q)\otimes\mb M_n(\mb C)$.  So if $(\iota\otimes{\rm id})(M)\geqslant 0$ for some $M\in\mc V\otimes\mb M_n(\mb C)$, then $M\in\left( \mc Q\otimes\mb M_n\right)^{\vee\vee}$, the double dual cone. Since $\mc Q\otimes \mb M_n$ is archimedean, this means $M+\epsilon \in \mc Q\otimes \mb M_n,$ and this implies $(\varrho\otimes{\rm id})(M)\geqslant 0.$ This proves that $\overline\varrho$ is u.c.p.
\end{proof}

Given $\varrho\colon \mc V\rightarrow \mb B(\mc H)$ $*$-linear, where $\mc H$ is of finite dimension $n$, we  define a functional \begin{align*} c_\varrho\colon \mc V\otimes \mb M_n(\mb C) &\rightarrow \mb C \\ v\otimes M&\mapsto  {\rm tr}(\varrho(v)M). \end{align*}

\begin{Prop}[Choi's Theorem] Let $(\mc A,\mc Q)$ be an archimedean quadratic module and $\mc V\subseteq\mc A$ a unital $*$-subspace.  Let $\mc H$ be a Hilbert space of dimension $n<\infty,$ and $\varrho\colon\mc V\rightarrow \mb B(\mc H)$ unital  $*$-linear. Then the following are equivalent:
\begin{itemize}
\item[(i)] $\varrho$ is u.c.p.
\item[(ii)] $\varrho\otimes{\rm id}$ maps $(\mc Q\otimes\mb M_n)_{\vert \mc V}$ to positive operators.
\item[(iii)] $c_\varrho$ is nonnegative on $(\mc Q\otimes\mb M_n)_{\vert\mc V}$.
\end{itemize}
\end{Prop}
\begin{proof}
Clear from the standard version of  Choi's Theorem \cite[Theorem 6.1]{pau} and the above considerations.
\end{proof}

We finally formulate the real closed separation theorem from \cite{netzthom} in the more general context of quadratic modules.

\begin{Def}
The quadratic module $(\mc A,\mc Q)$ is called {\it tame}, if $\mc Q=\bigcup_{i\in I} \mc Q_i,$ where \begin{itemize}
\item $(I,\leq)$ is a directed poset
\item each $\mc Q_i$ is a closed convex cone in a finite-dimensional subspace of $\mc A^h$
\item $i\leq j \Rightarrow \mc Q_i\subseteq \mc Q_j$ for all $i,j\in I$
\item for each finite-dimensional subspace $\mc V\subseteq \mc A$ and each $i\in I$ there exists $j\in I$ such that $\mc V^*\mc Q_i\mc V\subseteq \mc Q_j.$
\end{itemize}
\end{Def}

\begin{Ex}\label{tamex}
Assume $\mc Q\cap -\mc Q=\{0\}$ and $\mc Q$ admits a generating set $S\subseteq\mc Q$, such that $$v^*qv=0\Rightarrow v=0$$ holds for all $v\in\mc A, q\in S.$ Then $(\mc A,\mc Q)$ is tame. To see this, let $$I=\left\{ (\mc V, T)\mid \mc V \mbox{ finite-dimensional subspace of } \mc A, T\subseteq S \mbox{ finite}\right\}$$ be equipped with the obvious partial order. For $i=(\mc V,T)\in I$ we define $$\mc Q_i=\left\{\sum_{q\in T} \sum_{j} v_{qj}^*qv_{qj} \mid v_{qj}\in\mc V\right\}.$$ Using the arguments from \cite[Proposition 2.6 and Lemma 2.7]{posch}, closedness of $\mc Q_i$ follows if we show $$\sum_{q\in T}\sum_j v_{qj}^*qv_{qj}=0 \Rightarrow v_{qj}=0\quad  \forall q,j.$$ But this is clear from our assumptions.\end{Ex}

\begin{Thm}\label{rcsep}
Let $(\mc A,\mc Q)$ be a tame quadratic module, and $a\in \mc A^h\setminus \mc Q$. Then there exists an extension $\mb R\subseteq {\bf R}$ of real closed fields, a ${\bf C}$-vectorspace $\mc H$ with inner product, a  $*$-homomorphism of ${\bf C}$-algebras $$\pi\colon \mc A\otimes_{\mb C}{\bf C} \rightarrow {\mb L}(\mc H)$$ mapping $\mc Q\otimes{\bf R}_+$ to positive operators, and $\xi\in\mc H$ with  $$\langle \pi(a\otimes 1)\xi,\xi\rangle <0.$$
Furthermore, $\xi$ can be assumed to be cyclic w.r.t. $\pi$, meaning each $h\in\mc H$ is of the form $\pi(b)\xi$ for some $b\in \mc A\otimes{\bf C}.$
\end{Thm}
\begin{proof}
Since the argument is an adaption of the results in \cite{netzthom}, we skip the technical details.  For every $i\in I$ we find a linear functional $\varphi_i\colon \mc A^h\rightarrow \mb R$ with $\varphi_i\geq 0$ on $\mc Q_i$ and $\varphi(a)<0$. Choose an ultrafilter $\omega$ on $I$ containing all the the  upper sets $\{i\in I\mid i\geq j\}$, and let ${\bf R}=\mb R^\omega$ be the ultrapower. Then the $\mb R$-linear map $\varphi\colon\mc A^h\rightarrow{\bf R}; b \mapsto (\varphi_i(b))_{i\in I}$ separates $b$ from $\mc Q$ (using \L os's Theorem from model theory), and we can extend $\varphi$ to a unital $*$-linear map $\varphi\colon\mc A\rightarrow {\bf C}={\bf R}[i]$. Using the fourth property from the definition of a tame quadratic module (and \L os's Theorem again), one checks that $\varphi$ is u.c.p.  It follows that  the ${\bf C}$-linear map $\varphi\otimes{\rm id}\colon \mc A\otimes_{\mb C}{\bf C}\rightarrow {\bf C}$ is  nonnegative on $\mc Q\otimes {\bf R}_+$, and we perform the usual GNS construction, that works as well over $\bf C$. From this the result follows.
\end{proof}

\section{Residually finite-dimensional quadratic modules}\label{rfdsec}

We now define the notion of a residually finite-dimensional quadratic module, and characterize it in several ways. First, the notion of the Fell topology for $\ast$-representations of \Cs-algebras easily generalises to u.c.p.\ maps (see \cite[Section F.2]{be} for more details on the Fell topology).
\begin{Def}
 Let $\mc A$ be a \Cs-algebra and $\varrho\colon\mc A\to \mb B(\mc H)$ be a u.c.p.\ map. A functional of positive type associated with $\varrho$ is a functional
 \[
  \varphi(a)= \ip{\varrho(a)\xi,\xi},
 \]
where $\xi\in \mc H$ is a unit vector. We denote the set of such functionals by $\mathrm{Pos}(\varrho)\subset \mc A^*$.
\end{Def}

\begin{Def}
 Let $\mc A$ be a \Cs-algebra and $\varrho\colon \mc A\to \mb B(\mc H_\varrho)$ be a u.c.p.\ map, let $F\subset \mc A$ and $\Phi\subset \mathrm{Pos}(\varrho)$ be finite. The neighborhood system
 \begin{multline*}
  N(\varrho,F,\Phi,\eps)= \\\left\{\varrho'\colon \mc A\to \mb B(\mc H_{\varrho'})\text{ u.c.p.}\mid\forall \varphi\in\Phi  \ \exists\,\psi_i\in \mathrm{Pos}(\varrho')\ \forall a\in F: \left|\varphi(a) - \sum_{\mbox{ \tiny finite}} \psi_i(a)\right| < \eps \right\}
\end{multline*}
defines a topology on the set of u.c.p.\ maps $\mc A\to \mc B(\mc H)$ called the Fell topology.
We say that a u.c.p.\ map $\varrho'$ is weakly contained in $\varrho$, denoted $\varrho'\prec \varrho,$ if $\varrho'\in \overline{\{\varrho\}}$.
\end{Def}
The Fell topology just defined clearly coincides with the usual Fell topology when restricted to the set of representations of $\mc A$. The Stinespring Dilation Theorem implies at once that any u.c.p. map is weakly contained in its Stinespring dilation.

\begin{Def}
 An archimedean quadratic module $(\mathcal Q,\mathcal A)$ is residually finite-dimensional (r.f.d.) if its universal \Cs-algebra $C^*(\mathcal A,\mathcal Q)$ is residually finite dimensional, meaning that finite dimensional representations are dense in the set of all representations.\end{Def}

\begin{Thm}\label{pos}
For an archimedean quadratic module $(\mc A,\mathcal{Q})$, the following are equivalent:
\begin{itemize}

\item[(i)] $(\mc A,\mathcal{Q})$ is r.f.d.

\item[(ii)] For any $a\in \mathcal A^h$, whenever $\pi(a)\geqslant 0$ for all $\pi\in \Rep_{fd}(\mc A, \mathcal{Q})$, then $$a+\epsilon \in \mathcal{Q} \quad\forall \epsilon >0.$$
\end{itemize}
\end{Thm}
\begin{proof}
(i) $\Rightarrow$ (ii): if $C^*(\mc A,\mc Q)$ is r.f.d., $\pi(a)\geqslant 0$ for all $\pi\in \Rep_{fd}(\mc A, \mathcal{Q})$  is equivalent to $\pi(a) \geqslant 0$ for all $\pi\in\Rep(\mc A,\mc Q)$. By the abstract Positivstellensatz \cite{schm}, $a+\eps\in \mc Q$ for all $\eps > 0$ is equivalent to $\pi(a) \geqslant 0$ for all $\pi\in\Rep(\mc A,\mc Q)$.

(ii) $\Rightarrow$ (i): $a+\eps \in \mc Q$ for all $\eps > 0$ implies $\iota(a)\geqslant 0$ in $C^*(\mc A,\mc Q)$. If positivity in a \Cs-algebra is detected by finite-dimensional representations, then the \Cs-algebra is r.f.d.
\end{proof}

\begin{Rem}
Besides the above Positivstellensatz, property r.f.d.\ is also interesting from a computational point of view. The (semi)-norm $\Vert a\Vert_{\mc Q}$ of an element $a\in\mathcal A$ can be approximated in the following way (see also \cite{fnt} for more information): Upper bounds are obtained by computing numbers $\ell$ such that $\ell^2-a^*a\in\mc Q$. This is a semidefinite program, if only finitely many generators of $\mc Q$ and sums of squares from a finite-dimensional subspace of $\mc A$ are used. Making these constraints less and less restrictive, the sequence of upper bounds converges to $\Vert a\Vert_{\mc Q}$ from above.

Now a sequence of lower bounds is obtained by computing  $\sup \Vert \pi(a)\Vert$ over all representations $\pi$ of some bounded dimension. In case of a finitely generated quadratic module in a finitely generated algebra, this is a semialgebraic decision problem, which is decidable. If $\mc Q$ is r.f.d., these lower bounds will also converge to $\Vert a\Vert_{\mc Q}$, with growing dimension.
\end{Rem}

\begin{Lemma}\label{dense}
 A \Cs-algebra $\mc A$ is r.f.d.\ if and only if the set of its finite-dimensional representations is dense in the set of u.c.p.\ maps $\mc A\to\mb B(\mc H)$.
\end{Lemma}
\begin{proof}
 A \Cs-algebra $\mc A$ is r.f.d.\ iff if the set of its finite-dimensional representations is dense in the set of all representations, which is in turn dense in the set of all u.c.p.\ maps by the remark above.
\end{proof}

\begin{Def}
 A u.c.p.\ map $\varrho\colon \mc A\to\mb B(\mc H)$ is finite-dimensional if $\mc H$ is finite-dimensional. A u.c.p.\ map $\varrho\colon \mc A\to\mb B(\mc H)$ is strongly finite-dimensional if it possesses a finite-dimensional Stinespring dilation.
\end{Def}

\begin{Thm}\label{strdense}
For a unital \Cs-algebra $\mc A$, the following are equivalent:

\begin{itemize}\item[(i)] $\mc A$ is r.f.d.
\item[(ii)] For every finite-dimensional u.c.p.\ map  $\varrho\colon \mc A\to\mb B(\mc H)$, every finite $F\subset \mc A$ and every $\eps > 0$ there exists a strongly finite-dimensional u.c.p. map $\tilde\varrho\colon \mc A\to\mb B(\mc H)$ such that 
\[
\Vert \varrho(a)-\tilde\varrho(a)\Vert\leq \eps,\quad a\in F.
\]
\item[(iii)] Every state $\omega\colon \mc A\to\mb C$ is a weak$^*$-limit of states associated to finite-dimensional representations.
\end{itemize}
\end{Thm}
\begin{proof}
The equivalence of (i) and (iii) is well-known, and (ii) obviously implies (iii). The proof that (i) implies (ii) is essentially the argument from \cite[Proposition 2.2]{kech}.\end{proof}

\begin{Ex}
If $\mc A$ is commutative (with trivial involution), then every archimedean quadratic module in $\mc A$ is r.f.d. This follows for example via Theorem \ref{pos} from the commutative archimedean Positivstellensatz \cite[Theorem 5.4.4]{ma}, or can be deduced via functional calculus for commuting families of operators. 
\end{Ex}

\begin{Ex} Let $\Gamma$ be a group, and $\mc A=\mb C\Gamma$ the group algebra, equipped with $\mc Q=\Sigma^2\mc A^2.$ Then r.f.d.\ for $\mc Q$ is a well-studied property in group theory. There are groups which are r.f.d., for example free groups $\mb F_m$ \cite{choi} or surface groups \cite{lu}, and there are groups which are not r.f.d., like ${\rm SL}_n(\mb Z)$ for $n\geq 3$ (they have Kazhdan's Property (T) and thus each finite-dimensional representation is an isolated point in the Fell topology, see \cite{be}).
\end{Ex}

\begin{Ex}\label{klep} We consider the class of examples from \cite{he2}.
Let $\mc A=\mb C\langle z_1,\ldots, z_n\rangle$ be the free algebra with $z_i^*=z_i.$ Fix Hermitian matrices $M_1,\ldots,M_n\in \mb M_s(\mb C)^h$ and let $$\mc L= I_s+M_1z_1+\cdots + M_nz_n$$ be the associated linear matrix pencil. Then $$\mc Q=\left\{ \sum_j p_j^*p_j + q_j^*\mc L q_j\mid p_j,q_j\in \mc A^s\right\}$$ is a quadratic module in $\mc A$ which is r.f.d. This can be deduced from the Positivstellensatz in \cite{he2}, but we prove it directly, and in fact re-prove this Positivstellensatz with our method.  A representation $\pi\in\Rep(\mc A,\mc Q)$ is just a tuple $\ul T=(T_1,\ldots, T_n)$ of  self-adjoint operators on a Hilbert space $\mc H$, fulfilling $$\mc L(\ul T)=I_s\otimes {\rm id}_{\mc H}+ M_1\otimes T_1+\cdots + M_n\otimes T_n\geqslant 0$$ (to see this, use that $\pi$ can be assumed to admit a cyclic vector). For any isometry $\gamma\colon\tilde{\mc H}\rightarrow \mc H$ we set $\gamma^*\ul T\gamma=(\gamma^*T_1\gamma,\ldots,\gamma^*T_n\gamma)$ and find $$\mc L(\gamma^*\ul T \gamma)=(I_s\otimes\gamma)^*\mc L(\ul T)(I_s\otimes\gamma)\geqslant 0.$$ So the tuple $\gamma^*\ul T\gamma$ gives rise to a representation $\tilde\pi$  of $\mc Q$ again. The usual compression trick \cite[Theorem 6.1]{netzthom} shows that for finite-dimensional subspaces $\mc V\subseteq \mc A,$ $\mc H_0 \subseteq \mc H$, there is some isometry $\gamma\colon \mc H_1\rightarrow \mc H$ from a finite-dimensional space $\mc H_1\supseteq \mc H_0,$ such that \begin{align}\label{approx}p(\gamma^*\ul T\gamma)\equiv \gamma^* p(\ul T)\gamma\end{align} on $\mc H_0$,  for all $p\in\mc V$. Thus $\tilde\pi\in\Rep(\mc A,\mc Q)$ is close to $\pi$ in the Fell topology. This shows that $(\mc A,\mc Q)$ is r.f.d., even in a very strong sense.

We strengthen the argument to prove the strong Positivstellensatz from \cite{he2}. First check that $\mc Q$ is tame, using Example \ref{tamex}. Then for $a\in \mc A^h\setminus\mc Q$, use the real closed separation from Theorem \ref{rcsep}, and obtain $\pi(z_1),\ldots,\pi(z_n)\in\mb L(\mc H)$ with  $I_s\otimes {\rm id}_{\mc H}+ M_1\otimes \pi(z_1)+\cdots + M_n\otimes \pi(z_n)\geqslant 0$ (again use that $\pi$ admits a cyclic vector). Now apply the same compression trick as before to $\pi$, and obtain a finite-dimensional representation over ${\bf C}$ in which $a$ is not positive. Using Tarski's Transfer Principle, such a representation also exists over $\mb C$.  We have thus shown: If $a\in\mc A^h$ is nonnegative on $\Rep_{fd}(\mc A,\mc Q)$, then $a\in\mc Q$.
\end{Ex}

\begin{Ex}\label{putinar} This is the example from \cite{he5}.
 Let $\mc A=\mb C\langle z_1,z_1^*,\ldots, z_n,z_n^*\rangle$ and $$\mc Z=\left\{ (M_1,\ldots, M_n)\mid M_i \mbox{ matrices}, \sum_i M_i^*M_i=I\right\}.$$ Then $$\mc Q=\Sigma^2\mc A + \left\{ p\in\mc A^h\mid p\equiv 0 \mbox{ on } \mc Z\right\}$$ is r.f.d. This can be deduced from the Positivstellensatz in \cite{he5}, which we again re-prove it with our above separation method. Let $\pi\in\Rep(\mc A,\mc Q)$ be given. With $T_i=\pi(z_i)\in\mb B(\mc H)$ we have $\sum_i T_i^*T_i={\rm id}_{\mc H}.$ Let $\mc V\subseteq\mc A, \mc H_0\subseteq\mc H$ be finite-dimensional subspaces. Without loss of generality assume that $\mc V=\mc A_d$, the space of all polynomials of degree at most $d\geq 2.$ Inductively define $$\mc H_{i+1}={\rm span}\left\{ p(\ul T)h\mid p\in \mc V, h\in\mc H_i\right\}$$ for $i=0,1.$ Let $\gamma\colon \mc H_{2}\hookrightarrow\mc H$ be the embedding and consider the compressed operators $M_i:=\gamma^*T_i\gamma\in\mb B(\mc H_{2}).$ We have $p(\ul M)\equiv \gamma^*p(\ul T)\gamma$ on $\mc H_{1}$ for all $p\in\mc V$, and thus $$M\colon \mc H_1\rightarrow \mc H_{2}^n; h\mapsto (M_1h,\ldots,M_nh)$$ is an isometry. So we can extend to an isometry $$\tilde M=(\tilde M_1,\ldots,\tilde M_n)\colon \mc H_{2}\rightarrow \mc H_{2}^n,$$ and thus obtain a finite-dimensional representation $\tilde\pi$ of $(\mc A,\mc Q)$ on $\mc H_{2}.$ Now one checks that $p(\ul{\tilde M})\equiv \gamma^*p(\ul T)\gamma$ on $\mc H_0$, and so $\tilde \pi$ is close to $\pi$ in the Fell topology. So $(\mc A,\mc Q)$ is r.f.d.\ in a strong sense.

 Now we do the same over a real closed field $\bf R$. This time, we first pass to $$\mc B=\mc A/\left\{ p\in \mc A\mid p\equiv 0 \mbox{ on }\mc Z\right\},$$ where $\Sigma^2\mc B$ is tame, as is easily checked (using Example \ref{tamex}). We separate by a real-closed representation  and lift it to $\mc A$. Then we do the compression as described above, and transfer to $\mb C$ in the end. We have shown: If $a\in\mc A^h$ is nonnegative on $\Rep_{fd}(\mc A,\mc Q)$, then $a\in \mc Q$.
 
 \end{Ex}

 \begin{Ex}\label{more}
 Essentially the same methods can be used to show that the following quadratic modules are r.f.d., and even fulfill the strong Positivstellensatz $$a\in\mc A^h, a\geqslant 0 \mbox{ on } \Rep_{fd}(\mc A,\mc Q)\ \Rightarrow a\in \mc Q.$$
 \begin{itemize}
 \item[(i)] $\mc A=\mb C\langle z_1,z_1^*,\ldots, z_n,z_n^*\rangle$ with $\mc Q$ generated by either $$1-\sum_{i=1}^n z_i^*z_i\ \mbox{ or }\  1-z_i^*z_i \mbox{ for all } i.$$ Here we can separate and compress  without any further adjustments.
 \item[(ii)] $\mc A=\mb C\langle z_{ij}, z_{ij}^*\mid 1\leq i,j\leq n\rangle$ and  $ \mc Q=\Sigma^2\mc A+ \left\{ p\in\mc A^h\mid p\equiv 0 \mbox{ on } \mc Z\right\},$ where  $$\mc Z=\left\{ (M_{11},M_{12},\ldots, M_{nn})\mid M_{ij} \mbox{ matrices},  \left( M_{ij}\right)_{ij} \mbox{ unitary} \right\}.$$  After separating and compressing, we invoke Choi's matrix-trick from \cite{choi} and a suitable permutation of rows and columns.  
 \item[(iii)] $\mc A=\mb C\langle z_{ij},z_{ij}^*\mid 1\leq i,j\leq n\rangle$ and $\mc Q$ generated (as in Example \ref{klep}) by the quadratic matrix polynomial $$\mathcal P=I_n - (z_{ij})_{ij}^*(z_{ij})_{ij}\in \mb M_n(\mc A)^h.$$  This is even simpler as (ii).\end{itemize}
 \end{Ex}

 \begin{Ex}
 Let $\mc A=\mb C\langle u,u^*,v,v^*\rangle/(z^*z=zz^*=v^*v=vv^*=1)$ and let  $\mc Q$ be generated by $$\epsilon^2- (uv-vu)^*(uv-vu)$$ for some $\epsilon >0$ (this is called the {\it soft torus}). It is shown in \cite{eil} that $(\mc A,\mc Q)$ is rfd.
 \end{Ex}

 \begin{Ex}
 Let $\Gamma=\mb F_2\times\mb F_2$ and $\mc A=\mb C\Gamma$. Then $\mc A$ is rfd if and only if Connes' Embedding Conjecture is true (see \cite{kirch}).   
 \end{Ex}

\begin{Ex}\label{toep}
Let $\mc A=\mb C\langle z,z^*\rangle/(zz^*-1)$ be the Toeplitz algebra.  Then $\mc Q$ is not rfd. Finite dimensional representations correspond to unitary matrices, but the left-shift on $\ell^2(\mb N)$ yields a representation that cannot be approximated by finite-dimensional representations, since it is not unitary. 
\end{Ex} 

\section{Free Convexity}\label{convex}

Let us briefly introduce the main concepts of free convexity, as in \cite{he1,he3,he4,eff}. For some $n\geq 1$ we consider subsets $C_s\subseteq {\rm Her}_s(\mb C)^n$ of $n$-tuples of hermitian matrices of size $s$, for all  $s\geq 1.$ The whole collection $C=\bigcup_{s\geq 1} C_s$ is called matrix-convex, if it is closed under block-diagonal sums and compressions via isometries. That means, whenever $\ul A=(A_1,\ldots,A_n)\in C_s, \ul B=(B_1,\ldots,B_n)\in C_r$, and $V\in \mb M_{s,r}(\mb C)$ with $V^*V=1,$ then $$\ul A\oplus\ul B=\left( \left(\begin{array}{cc}A_1 & 0 \\0 & B_1\end{array}\right),\ldots, \left(\begin{array}{cc}A_n & 0 \\0 & B_n\end{array}\right)\right)\in C_{s+r}$$ and $$V^*\ul AV=\left( V^*A_1V,\ldots, V^*A_nV\right)\in C_r.$$
This easily implies that each $C_s$ is convex in the real vectorspace ${\rm Her}_s(\mb C)^n$, but matrix convexity of $C$ is a stronger assumption in general.

For any set  $C=\bigcup_s C_s,$ its matrix-convex hull ${\rm mconv}(C)$ is the smallest matrix-convex superset of $C$. In case that $C$ is already closed under block-diagonal sums, it is easy to see that we only need to add compressions to obtain the matrix convex hull: \begin{align}\label{comp}{\rm mconv}(C)_s=\left\{ V^*\ul AV\mid r\geq s, \ul A\in C_r, V\in\mb M_{r,s}(\mb C), V^*V=1 \right\}.\end{align} 
Now assume $\mc A=\mb C\langle z_1,\ldots, z_n\rangle$ with $z_i^*=z_i$ and $p\in \mc A^h$. Define $$C(p)_s=\left\{ (A_1,\ldots,A_n)\in {\rm Her}_s(\mb C)^n\mid p(A_1,\ldots,A_n)\geqslant 0\right\}$$and  $C(p)=\bigcup_s C(p)_s$, a so-called (basic closed)  free semialgebraic set (finite intersections of such sets are also called basic closed). 
Understanding the matrix convex hullof such sets is one of the main issues in the above mentioned papers. 

Note that one can also use operators instead of matrices to define free convex hulls. For a Hilbert space $\mc H$ define $$C(p)_{\mc H}=\left\{ (T_1,\ldots, T_n)\mid T_i \in\mb B(\mc H)^h, p(T_1,\ldots,T_n)\geqslant 0\right\}$$ and the operator convex hull ${\rm oconv}(C(p))$ as $${\rm oconv}(C(p))_s=\left\{ V^* \ul T V \mid \mc H \mbox{ Hilbert space }, \ul T\in C(p)_{\mc H}, V\colon \mb C^s\rightarrow \mc H \mbox{ isometry}\right\}.$$ 
Now, interesting (and previously open) questions are:
\begin{itemize}
\item Can $r$ in (\ref{comp}) be bounded in terms of $s$ (and maybe other data)?
\item Is ${\rm mconv}(C(p))$ and/or ${\rm oconv}(C(p)) $ semialgebraic in any (free) sense?
\item Is at least each ${\rm mconv}(C(p))_s$ and/or ${\rm oconv}(C(p))_s$ semialgebraic in the usual sense?
\end{itemize}

We will answer these questions to the negative below. But let us first define a broader and coordinate-free notion of free convexity. Example \ref{toe} will show that this might be useful. Note that all of the following concepts coincide for an archimedean quadratic module  $(\mc A,\mc Q)$ and its universal $C^*$-algebra $C^*(\mc A,\mc Q).$ 

\begin{Def} Let $\mc A$ be a $C^*$-algebra.
The convex hull of $\Rep(\mc A)$ is defined as $${\rm conv}\Rep(\mc A)=\left\{ \varrho\colon \mc A\rightarrow\mb B(\mc H)\mid \mc H \mbox{ Hilbert space}, \varrho \mbox{ u.c.p.}\right\}.$$
 The convex hull of $\Rep_{fd}(\mc A)$ is defined as  $${\rm conv}\Rep_{fd}(\mc A)=\left\{ \varrho\colon \mc A\rightarrow\mb B(\mc H)\mid \varrho \mbox{ strongly finite dimensional u.c.p.}\right\}.$$
\end{Def}

\begin{Rem}
(i) Convex hulls are the smallest supersets closed under compressions via isometries. This immediately follows from Stinespring's Dilation Theorem. Also note that  the sets $\Rep(\mc A), \Rep_{fd}(\mc A)$ and their convex hulls are closed under finite direct sums.

(ii) In case $\mc A=\mb C\langle z_1,\ldots, z_n\rangle$, we obtain the old notions of free convex hulls when we restrict the u.c.p.\ maps to the space $\mc V={\rm span}\{z_1,\ldots, z_n\}.$
\end{Rem}

 \begin{Ex}\label{toe}
 Let $\mc A=\mb C\langle z,z^*\rangle/(zz^*-1)$ be the Toeplitz algebra with $\mc Q=\Sigma^2\mc A,$ and $\mc V={\rm span}\{ z,z^*\}.$ For any finite-dimensional $\varrho\in {\rm conv}\Rep(\mc A,\mc Q)$ there is a strongly finite dimensional $\tilde\varrho\in {\rm conv}\Rep_{fd}(\mc A,\mc Q)$ with $\varrho\equiv\tilde\varrho$ on $\mc V$. In fact $\varrho(z)=\gamma^*\pi(z)\gamma$ for some $\pi\in\Rep(\mc A,\mc Q)$ and some isometry $\gamma.$ So $\varrho(z)$ is a  finite-dimensional contration, and thus admits a finite-dimensional unitary dilation.

 On the other hand, $\mc Q$ is not rfd, as shown in Example \ref{toep}, so for other subspaces $\mc V$ we don't even get a good approximation by strongly finite-dimensional morphisms (by Theorem \ref{strdense}).  This suggests that restricting all the maps to a generating subspace $\mc V$ of $\mc A$, as done in free convexity throughout so far, is not always a good idea. This approach will for example not be able to detect the property rfd. This is why we proposed  the above notion of free convex hulls.

 \end{Ex}

Given a set $\mc T$ of mappings defined on $\mc A$,  and a subset $\mc V\subseteq \mc A$, we call $$\mc T_{\mid\mc V}=\left\{ \varrho_{\mid\mc V}\mid \varrho\in\mc T\right\}$$ the projection of $\mc T$ to $\mc V$. Note that in a finite-dimensional real vectorspace, there is a notion of {\it semialgebraic set}, which is independent of the choice of a basis.

\begin{Thm}\label{projsa}
Let $\mc A$ be a $C^*$-algebra, $\mc V\subseteq\mc A$ a finite-dimensional unital $*$-subspace and $\mc H$ a Hilbert space with $\dim(\mc H)=n<\infty$. Then the following are equivalent:
\begin{itemize}
\item[(i)] The projection of $\left\{ \varrho\colon\mc A\rightarrow \mb B(\mc H)\mid \varrho \in {\rm conv}\Rep(\mc A)\right\}$ to $\mc V$ is semialgebraic.
\item[(ii)] $(\Sigma^2\mc A \otimes \mb M_n)_{\mid \mc V}$ is semialgebraic.
\end{itemize}
\end{Thm}
\begin{proof}
A unital $*$-linear mapping $\varrho\colon \mc V\rightarrow \mb B(\mc H)$ is in the projection from (i), if and only if it is u.c.p, by Arvesons's Extension Theorem. By Choi's Theorem, this is equivalent to $c_\varrho$ being nonnegative on $(\Sigma^2\mc A \otimes \mb M_n)_{\mid\mc V}.$ So the projection from (i) is just the dual of the closed set $(\Sigma^2\mc A\otimes \mb M_n)_{\mid \mc V}$. This proves the claim.
\end{proof}

 \begin{Ex} Let $\mc A$ be a commutative $C^*$-algebra. An element from $\mc V\otimes \mb M_n(\mb C)$ is positive  if and only if it is positive under each representation of the form $\pi\otimes{\id}$, where $\pi\in\Rep(\mc A)$. Since $\mc A$ is rfd we can restrict to $\pi\in\Rep_{fd}(\mc A)$, and since $\mc A$ is commutative, even to one-dimensional representations.
 
If $\mc A$ is the universal $C^*$-algebra of a finitely generated quadratic module in a finitely generated commutative algebra, then the set of one-dimensional representations is semialgebraic. Thus $(\Sigma^2\mc A\otimes \mb M_n)_{\mid \mc V}$ and the corresponding projection of ${\rm conv}\Rep(\mc A)$ to $\mc V$ are always semialgebraic.
 \end{Ex}

 \begin{Ex}
 Let $\mc A=\mb C\mb F_m$ be the group algebra of the free group. Again $(\Sigma^2\mc A \otimes \mb M_n)_{\mid\mc V}$ is always semialgebraic. As before, an element from $\mc V \otimes \mb M_n(\mb C)$ is positive if and only if it is positive at each representation $\pi\otimes{\rm id},$ where $\pi\in\Rep(\mc A).$ Choi's proof that $\mc A$ is rfd shows that we can even restrict to representations $\pi$ of some fixed dimension, depending only on $\mc V$ and $n$. So $(\Sigma^2 \mc A\otimes \mb M_n)_{\mid \mc V}$ can be defined by a formula in the language of ordered rings and is thus semialgebraic.
 The same reasoning applies to all quadratic modules from Examples \ref{klep}, \ref{putinar} and \ref{more}.
 \end{Ex}
 
 In general, the projections of both ${\rm conv}\Rep(\mc A)$ and ${\rm conv}\Rep_{fd}(\mc A)$ are not semialgebraic, answering the above questions.
 
 \begin{Thm} There exists a finitely generated quadratic module $\mc Q$ in the free algebra $\mc A=\mb C\langle z_1,\ldots, z_n\rangle$, such that already the projection of $$\left\{ \rho\colon \mc A\rightarrow \mb C\mid \rho\in {\rm conv}\Rep_{fd}(\mc A,\mc Q)\right\}$$
to $\mc V={\rm span}\{z_1,\ldots, z_n\}$ is not semialgebraic.
 \end{Thm}
 \begin{proof} Instead of the free algebra, we work with the group algebra $\mc A=\mb C\Gamma$ of the discrete Heisenberg group $\Gamma=\langle a,b,c\mid c=aba^{-1}b^{-1}, ca=ac, cb=bc\rangle$ and $\mc Q=\Sigma ^2\mc A$. By lifting the relations as pairs of inequalities to the free algebra, one obtains an example in the free algebra. 
 
 Each irreducible $n$-dimensional representation of $\mc A$ maps $c$ to an $n$-th root of unity. This is true since $c$ lies in the center of $\mc A$, and as a commutator has determinant one. Any $n$-th root of unity is attained through a representations, by \cite{fan}. Let $\tilde{\mc V}={\rm span}\{(c+c^*)/2,(c-c^*)/(2i)\}\subseteq \mc A$. Then the projection of $\left\{ \varrho\colon\mc A\rightarrow \mb C\mid \varrho\in {\rm conv}\Rep_{fd}(\mc A)\right\}$  to $\tilde{\mc V}$ is$${\rm conv}\left\{ (x,y)\in\mb R^2\mid x+iy \mbox{ roof of unity}\right\},$$ which is not semialgebraic.\end{proof}
 
 \begin{Rem} The  example also shows that there is no bound on $r$ in (\ref{comp}). 
 \end{Rem}
 
\begin{figure}[h]

 \includegraphics[scale=0.3]{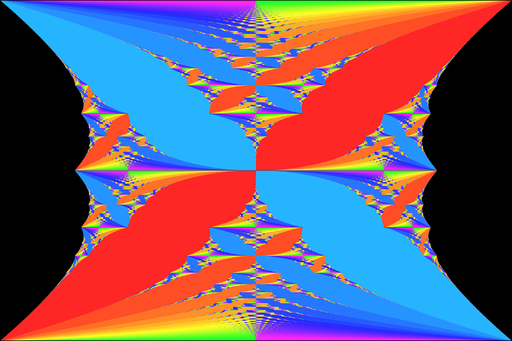}
  \caption{\label{butt} The Hofstadter butterfly \newline \small\url{https://commons.wikimedia.org/wiki/File\%3AHofstadter's\_butterfly.png}}
 \end{figure}

\begin{Thm}\label{heis2}  There exists a finitely generated quadratic module $\mc Q$ in the free algebra $\mc A=\mb C\langle z_1,\ldots, z_n\rangle$, such that already the projection of $$\left\{ \rho\colon \mc A\rightarrow \mb C\mid \rho\in {\rm conv}\Rep(\mc A,\mc Q)\right\}$$
to $\mc V={\rm span}\{z_1,\ldots, z_n\}$ is not semialgebraic.

\end{Thm}
\begin{proof}
Again we work in the group algebra $\mb C\Gamma$ of the discrete Heisenberg group $\Gamma=\langle a,b,c\mid c=aba^{-1}b^{-1}, ca=ac, cb=bc\rangle$. This time let $\tilde{\mc V} = \mathop\mathrm{span}\{a+a^*+b+b^*,\,(c+c^*)/2,\,(c-c^*)/(2i)\}$. The classification of irreducible representations of $\Gamma$ is well-known, the spectral properties of the above operators in these representations are extensively studied in \cite{bvz}, and we will use these results. The irreducible representations of $\Gamma$ are parametrised by the circle $\{e^{i\theta}\mid \theta\in[0,2\pi]\}$, and in every such irreducible representation we have $\pi_\theta(c) = e^{i\theta}$, $\pi_\theta(ab) = e^{i\theta} \pi_\theta(ba)$, so in an irreducible representation $a$ and $b$ generate a noncommutative torus with parameter $\theta$. We denote $H_\theta = \pi_\theta(a + a^* + b + b^*)$. Using the automorphism of the noncommutative torus which maps the generators to the negatives of them, it is not hard to see that the spectrum of $H_\theta$ is symmetric.
 
Now a u.c.p. map $ \varrho\colon\mc A\rightarrow \mb C$ is just a state on $C^*(\Gamma)$, and it's a well-known general fact that the states on a \Cs-algebra form a closed convex set whose extremal points are the pure states coming from irreducible representations. Thus, the projection of $\left\{ \varrho\colon\mc A\rightarrow \mb C\mid \varrho\in {\rm conv}\Rep(\mc A)\right\}$ to $\tilde{\mc V}$ is the closed convex hull
 \[
  C = \overline{\rm conv}\left\{ (\pm\norm{H_\theta},\cos\theta,\sin\theta)\in\mb R^3\mid \theta\in[0,2\pi]\right\}.
 \]

The function $\theta \mapsto \norm{H_\theta}$ describes the boundary of the ``Hofstadter butterfly'' \cite{ho} (see Figure \ref{butt} for a picture), and is known to be non-differentiable at the points where $\theta/(2\pi)$ is rational \cite{ram,hs}.  
So if $C$ were a semialgebraic set, its intersection with the cylinder $Z = \{(x,y,z)\in\mb R^3\mid y^2 + z^2 = 1\}$ would also be semialgebraic, and thus the functions $\theta \mapsto \pm\norm{H_\theta}$ whose graphs form the (relative) boundary of the set $Z\cap C$ would be piecewise smooth, which yields a contradiction. Therefore $C$ is not semialgebraic.
\end{proof}

 The membership problem from real algebraic geometry is the following: Given a finitely generated quadratic module $(\mc A,\mc Q)$ and a finite-dimensional $\mb R$-subspace $\mc V\subseteq \mc A^h$, is $\mc Q\cap \mc V$ a semialgebraic set? This is known to be true in certain cases, but an open question in general \cite{aug}. 
 
 \begin{Cor} There is a finitely generated quadratic module in the free algebra $\mb C\langle z_1,\ldots, z_n\rangle,$ for which the membership problem has a negative answer.
 \end{Cor}
 \begin{proof} This follows from
 Theorem \ref{heis2}  combined with Theorem \ref{projsa}.   In fact, if $\mc Q\cap\mc V$ is semialgebraic, then so is its closure $\overline{\mc Q\cap\mc V}$, and this is equivalent to condition (ii) for the universal $C^*$-algebra $C^*(\mc A,\mc Q)$ in Theorem \ref{projsa}.
 \end{proof}
 
\section*{Acknowledgments}
This research was supported by ERC Starting Grant No.\ 277728.

\end{document}